\input amstex
\documentstyle{amsppt}
\def\pf{\hfill\hfill\qed}

\def\br{\Bbb R}

\leftheadtext{Augustin Banyaga}
\TagsOnRight
\NoBlackBoxes
\topmatter
\title
A HOFER-LIKE METRIC ON THE GROUP OF SYMPLECTIC DIFFEOMORPHISMS  
\endtitle
\vskip .3in
\author
Augustin Banyaga
\endauthor
\vskip .3in
\abstract
Using a "Hodge decomposition" of symplectic isotopies on a compact symplectic manifold $(M,\omega)$
, we construct a norm on the identity component in the group of all
symplectic diffeomorphisms of $(M,\omega)$ whose restriction to the group $Ham(M,\omega)$ of
 hamiltonian diffeomorphisms is  bounded from above by the Hofer norm. Moreover, 
 $Ham(M,\omega)$  is closed in $Symp(M,\omega)$ equipped with the topology 
induced by the extended norm. We give an application to the $C^0$ symplectic topology.
 We also discuss extensions of Oh's spectral distance.

\endabstract
\endtopmatter
\document
\baselineskip 20pt

{\bf 1. Introduction and statement of the main results}

\vskip .1in

Let $Symp(M,\omega)$  denote the group of  all symplectic diffeomorphisms of a
 compact symplectic manifold $(M,\omega)$, endowed with the $C^{\infty}$ compact-open topology,
 and $Symp(M,\omega)_0 = G_{\omega}(M)$ the identity component in $Symp(M,\omega)$. $Symp(M,\omega)_0$
 consists of  symplectic diffeomorphisms $h$ such that there is a symplectic isotopy
 $h_t$ from the identity to $h$. By definition $h_t$ is a symplectic isotopy if the 
 map $(x,t) \mapsto h_t(x)$ is smooth and for all t , $h_t^*\omega = \omega$. We denote by $Iso(M)$ the set 
 of all symplectic isotopies, and by $Iso(\phi)$ the set of all symplectic isotopies
 from the identity to $\phi \in Symp(M,\omega)_0$.
\vskip .1in 
  Let $Ham (M,\omega) \subset Symp(M,\omega)_0$ 
 be the subgroup of 
 Hamiltonian diffeomorphisms. A  diffeomorphism $\psi$ is Hamiltonian iff it is the time 1 map of  a smooth family 
 of diffeomorphisms $\psi_t$ such that if 
 $$
\dot {\psi}_t(x) = \frac{d\psi_t}{dt}(\psi^{-1}_t(x)),~~~~~ \psi_0(x) = x~~~~~~~~~~~~~~~~~~~~~~~~~~~~~~~~~~ \tag 1 
$$
there exists a smooth family of functions $u_t$ such that
$$
i_{(\dot {\psi}_t)} \omega =du_t.  ~~~~~~~~~~~~~~~~~~~~~~~~~~~~~~~~~~~~~~~~~~~~~\tag 2
$$
The family of diffeomorphisms $\psi_t$ above is called a hamiltonian isotopy.

\noindent We denote by $HIso(\phi)$ the set of all hamiltonian isotopies
 from $\phi \in Ham(M,\omega)$ to the identity, and by $HIso(M)$ the set of all hamiltonian isotopies.

\vskip .1in
                   
\noindent In  equation (2), $i_{(.)}$ denotes the interior product: $i_X\omega$ is the 1-form such 
that $i_X\omega(Y) = \omega(X,Y)$ . Recall that a symplectic form is a closed 2-form $\omega$ such that 
the map assigning to a vector field $X$ the 1-form $i_X\omega$ is an isomorphism $\tilde {\omega}$. For any 1-form 
$\alpha$, we denote by $\alpha^\#$ the vector field $\tilde {\omega}^{-1}(\alpha)$.

\vskip .1in

The {\it Hofer length} of a hamiltonian isotopy $\psi_t$ is defined as:
$$
l_H(\psi_t) =\int_0^1  ( max_{x} u_t(x) - min _{x} u_t(x)) dt ~~~~~~~~~~~~~~~~~~~~~~~~~\tag 3
$$

\vskip .1in

One also denotes 
$$
 max_{x} u_t(x) - min _{x} u_t(x)) = osc (u_t(x))
$$
and call it the oscillation of $u_t$.
 
 \vskip .1in
 
Hence the Hofer length is the mean oscillation of the hamiltonian $u_t$ 
of the hamiltonian isotopy $\Phi = (\phi_t)$.

\vskip .1in

For $\psi \in Ham(M,\omega)$, the {\it Hofer norm} is defined as: 
$$
||\psi||_H = inf (l_H(\psi_t))~~~~~~~~~~~~~~~~~~~~~~\tag 4
$$
where the infimum is taken over all hamiltonian isotopies $\psi_t\in HIso(\psi)$
and $u_t$ is the function in equation (2).
\vskip .1in

The  Hofer distance between two hamiltonian diffeomorphisms $\phi$ and $\psi$ is:
$$
d_H (\phi,\psi) = || \phi\psi^{-1}||_H
$$
It is easy to see that the formula above defines a bi-invariant pseudo-metric 
 but it is very challenging to show that it is  not degenerate and hence it is  
a genuine distance [5],[7],[12], [13]. 
\vskip .2in

In this paper we  propose a formula  for the length of a symplectic isotopy $\Phi = (\phi_t)$ (5), 
which generalizes  the length of a hamiltonian isotopy (3).

\vskip .1in

Fix a riemannian metric on $M$ and consider the Hodge decomposition of $i_{(\dot {\phi}_t)}\omega$ 
$$
 i_{(\dot {\phi}_t)}\omega = \Cal H^{\Phi}_t + du^{\Phi}_t  
$$
 
where  $\Cal H^{\Phi}_t$ and $u^{\Phi}_t$ are smooth family of harmonic 1-forms and functions respectively.
\vskip .1in

We define the length $l(\Phi)$ of the isotopy $\Phi$ by:
 $$
 l(\Phi) = \int_0^1 (|\Cal H^{\Phi}_t | + (max_x (u^{\Phi}_t) - min_x (u^{\Phi}_t))dt ~~~~~~~~~~\tag 5
 $$
 
 Here $|\Cal H^{\Phi}_t |$ is the "Euclidean" norm of the harmonic 1-form $\Cal H^{\Phi}_t $
  ( see (13), (14)).
  
 \vskip .1in
 
 This formula reduces to (3) for hamiltonian isotopies.
Unfortunately, unlike (3), we do not have
 $$
l(\Phi) = l(\Phi ^{-1})
$$
where $\Phi ^{-1} = (\phi_t^{-1})$. We will also write:
$$
 l(\Phi)
 = ||\dot {\phi}_t)||
$$
 \vskip .1in
 For any $\phi \in Symp(M,\omega)_0$, 
 we define the energy $e_0(\phi)$ of $\phi$ as:
 
 $$
 e_0(\phi) = inf_{\Phi \in Iso(\phi)}  (l(\Phi))
 $$
 
 \vskip .1in
 
 Our main result is the following
 
 \vskip .1in
 
 \proclaim{Theorem 1}
 
 Let $(M,\omega)$ be a closed symplectic manifold. Consider the map
  $e: Symp(M,\omega)_0 \to \br \cup \{\infty\}$
 $$
 e(\phi) = 1/2(e_0(\phi) + e_0(\phi^{-1})).
 $$
 
 Then $e$ is a norm on $Symp(M,\omega)_0$ whose restriction  to $Ham(M,\omega)$ is 
 bounded from above by the Hofer metric.

Moreover the subgroup $Ham(M,\omega)$ is closed in 
$Symp(M,\omega)$ endowed with the metric topology defined by $e$.

\endproclaim

\vskip .1in

We define a distance on $Symp(M,\omega)$ by:
$$
d(\phi, \psi) = e(\phi \psi^{-1})
$$
This distance is obviously right invariant, but not left invariant.

\vskip .1in

{\bf Remark}

 The fact that (5) reduces to (3) when $\Phi $ is a hamiltonian isotopy implies that
 $$
 e(\phi) \leq ||\phi||_H
 $$
 for all $\phi \in Ham(M,\omega)$.
 
 \vskip .1in

{\bf Conjecture} The restriction of the norm $e$ to $Ham(M,\omega)$ is equivalent to the Hofer norm.

\vskip .1in

What is the interest of our construction? From the Hofer norm, there are easy ways of constructing
bi-invariant norms on $Symp(M,\omega)$. One is given by Han [4]:

fix a positive number $K$ and define
$$
||\phi||_K = \cases min(||\phi||_H , K), &\text{if $\phi \in Ham(M,\omega)$}\\
K &\text{otherwise}.\endcases
$$

Another is given by Lalonde-Polterovich [8]:

fix a real number $\alpha$ and define
$$
||\phi||_{\alpha} = sup\{ ||\phi f \phi^{-1} f^{-1}||_H  | f \in Ham(M,\omega), ||f||_H \leq \alpha\}.
$$

In both cases the restriction of these metrics back to $Ham(M,\omega)$ gives different topologies
 on $Ham(M, \omega)$. In particular $Ham(M,\omega)$ in these topology has always a finite diameter
 which is known to be untrue for the Hofer norm in several cases.
 
Hence the advantage of our construction is that its restriction to $Ham(M,\omega)$  gives a "better" topology,
which may  be the same   if the conjecture is true.

\noindent Moreover the " Hofer-like" formula (5) allows to define a distance $D$
 on the space $Iso(M)$ of symplectic isotopies of $(M,\omega)$. If $\Phi = (\phi_t)$ and $\Psi = (\psi_t)$
 are symplectic isotopies define:
 $$
 D(\Phi, \Psi) = ||\dot {\phi}_t - \dot {\psi}_t|| = :\int_0^1 (| \Cal H^{\Phi_t} -\Cal H^{\Psi_t}| + osc
  (u^{\Phi_t} - u^{\Psi_t})) dt. ~~~~~~~~~~~~~~~~~~~~~~~\tag 5'
 $$
 
 It is clear that $D$ is a distance. Moreover if $\Phi, \Psi$ are hamiltonian isotopies, then
 $$
 D(\Phi, \Psi) = l_H(\Phi\Psi^{-1}).
 $$
 In general the formula above is not true.
 
\noindent The distance formula will be used to the define
 the {\it symplectic topology} on $Iso(M)$.

\vskip .1in

{\bf 2. Hamiltonian and harmonic  diffeomorphisms}

\vskip .1in

For each symplectic isotopy $\Phi = (\phi_t)$, consider the following 1-form:
$$
\Sigma(\Phi) = \int_o^1 (i_{\dot {\phi}_t}\omega) dt ~~~~~~~~~~~~~~~~~~~~~\tag 6
$$
It is shown  in [1], ( see also [2]) that the cohomology class $[\Sigma(\Phi] \in H^1(M,\br)$ 
of the form $\Sigma(\Phi)$ depends only on the class $[\Phi]$ of $\Phi$ in the universal covering $\tilde
 {G}(M,\omega)$ of $Symp(M,\omega)_0 = G(M,\omega)$ and that
the map $[\Phi] \mapsto [\Sigma(\Phi)]$ is a surjective homomorphism
$$
\tilde S :\tilde {G}(M,\omega) \to H^1(M,\br) ~~~~~~~~~~~~~~~~~~~~~~~~~\tag 7
$$
The group 
$$
\Gamma = \tilde S(\pi_1(G(M,\omega))) \subset H^1(M,\br)
$$
is called the {\bf flux group}.

\vskip .1in
In [1], it was observed that $\Gamma$ was discrete in several examples and the author wrote " I do 
not know any flux group which is not discrete". The statement that " $\Gamma$ is discrete" became known as the
 "Flux conjecture". This conjecture has been recently proved by Ono [11] using Floer-Novikov homology.
 
 \vskip .1in
 
 \proclaim
 {Theorem (Ono)}
 
 Let $(M,\omega)$ be a compact symplectic manifold, then the flux group is discrete.
 
 \endproclaim
 
 Consider the induced homomorphism :
 $$
 S : G(M,\omega) \to H^1(M,\br)/\Gamma ~~~~~~~~~~~~~~~~~~~~~~\tag 8
 $$
 
 \vskip .1in
 
 In [1], [2], it is shown that the Kernel of S coincides with the group $Ham(M,\omega)$ 
 of Hamiltonian diffeomorphisms, and it is a simple group, which coincides with the commutator subgroup
 $[G(M,\omega),G(M,\omega)]$ of $G(M,\omega)$. We summarize:
 $$
 Ham(M,\omega) = Ker S = [G(M,\omega),G(M,\omega)]~~~~~~~~~~~~~~~~~~~~~~\tag 9
 $$
 for all closed symplectic manifolds $(M,\omega)$.

 \vskip .1in
 
 We will need to represent in a unique way cohomology classes ; this is achieved by Hodge theory on 
 compact riemannian manifolds.
 
The Hodge decomposition theorem ( see for instance [14]) asserts that 
any smooth family of p-forms $\theta_t$ on a compact oriented riemannian manifold $(M, g)$ 
can be decomposed in a unique way as
$$
\theta_t = \Cal H_t + d\alpha_t + \delta \beta_t ~~~~~~~~~~~~~~~~~~~~~\tag 10
$$
where $\Cal H_t$ is harmonic, i.e $d\Cal H_t = \delta \Cal H_t = 0$ .
Here $\delta$ denotes the codifferential.

\noindent If $d\theta_t = 0$ , then $\delta \beta_t = 0$ .
\noindent The forms $\Cal H_t$, $\alpha_t$ and $\beta_t$ depend smoothly on $t$.

The harmonic form $\Cal H_t$ is  a unique representative of
the cohomology class $[\theta_t]\in H^1(M,\br)$ of $\theta_t$.

\vskip .1in

\proclaim
{Definition 1}

Let $(M,\omega)$ be a  compact symplectic manifold, equipped with some riemannian metric.
A vector field $X$ on $M$ is said to be 
a harmonic vector field if $i_X\omega$ is a harmonic form.

 A diffeomorphism $\phi$ of $M$ is said to be  a harmonic diffeomorphism  if
if there exists a smooth family $\Cal H_t$ of harmonic 1-forms
 such that $\phi$ is the time 1 map of the symplectic isotopy $\phi_t$ such that
 $$
 \dot {\phi}_t  = (\Cal H_t)^\#. ~~~~~~~~~~~~~~~~~~~~~~~~~~~~~~\tag 11
 $$
 We say that $\phi_t$ is a harmonic isotopy. 
 
 \endproclaim
 
 \vskip .1in 
 
 Let $symp(M,\omega)$ be the set of symplectic vector fields,
  $harm(M)$ the set of harmonic vector fields and $ham(M,\omega)$ the space of 
  hamiltonian vector fields. If $X \in symp(M,\omega)$
  then $i_X\omega$ is closed. The decomposition $i_X\omega = \Cal H + du$
  expresses  $X$ as
  $$
  X = H + X_u ~~~~~~~~~~~~~~~~~~~~~~~~~~~~\tag 12
  $$
  where $H = (\Cal H)^\#$ is harmonic and $X_u$ is the hamiltonian vector 
  field with hamiltonian $u$.
  
  \vskip .1in
  
 Hence $symp(M,\omega)$ is the Cartesian product of $harm(M)$ and 
 $ham(M,\omega)$. We give $symp(M,\omega)$  the product metric :
 $$
 |X| = |H| + max_x u(x) - min_x u(x) ~~~~~~~~~~~~~~~~~~~~~~~~\tag 13
 $$
 where $|H|$ is the norm on $harm(M)$ given below:
 
 \vskip .1in
 
 the space $harm(M)$, which is isomorphic to the space of harmonic 1-forms 
 is a finite dimensional vector space whose dimension is
 the first Betti number of $M$.
 
 \vskip .1in
 
 In this paper , we fix a basis $h_1,..., h_r$ of
  harmonic 1-forms and consider
 $(H_i) = (h_i ^{\#})$ the corresponding basis of $harm(M)$. We give these 2 vector spaces the following
  Euclidean metric : if $h = \sum_i \lambda_i h_i$, $H =  \sum_i \lambda_i H_i$,  then 
  $$
   |h| =  |H| =\sum_i |\lambda_i| ~~~~~~~~~~~~~~~~~~~~~~~~~~~~~~~\tag 14
 $$
  
 \vskip .1in
 
 In view of (13), the length formula (5) gives a Finsler metric on $Symp(M,\omega)$.
 
 \vskip .1in
 
 {\bf Remark}
 
 \vskip .1in
 
 The function $u$ in the Hodge decomposition $ i_X\omega = \Cal H + du$ is not necessarly normalized.
 However if in (13) $|X| = 0$, then $|H| = 0$, i.e  $i_X\omega = du$ and  $osc(u) =0$ 
 implies that $u$ is constant, and
  hence $du = 0$. Therefore $X = 0$.

\vskip .1in
 
 \proclaim
 {Lemma 1}
 
Any symplectic isotopy $\Phi = (\phi_t)$ on a compact symplectic manifold $(M,\omega)$ can be decomposed in a unique way as
 $$
 \phi_t = \rho_t.\psi_t
 $$
 where $\rho_t$ is a harmonic isotopy and $\psi_t$ is a hamiltonian isotopy. In particular,
  if $\phi_t$ is a hamiltonian isotopy, then $\phi_t = \psi_t$ and $\rho_t = id_M$.
\endproclaim 

 \vskip .1in
 
 \demo
 {Proof}
 
 By Hodge decomposition theorem $i_{(\dot {\phi}_t)}\omega$ can be decomposed in a unique way as
 $$
 i_{(\dot {\phi}_t)}\omega = \Cal H^{\Phi}_t + du^{\Phi}_t  
 $$
 
where  $\Cal H^{\Phi}_t$ and $u^{\Phi}_t$ are smooth family of harmonic 1-forms and functions respectively.
Let $\rho_t$ be the harmonic isotopy such that $\dot {\rho}_t = (\Cal H_t)^\#$.
Set now $\psi_t = (\rho_t)^{-1}.\phi_t$.
 From $\phi_t = \rho_t.\psi_t$ , we get: 
 $$
 \dot {\phi}_t = \dot {\rho}_t + (\rho_t)_* \dot {\psi}_t ~~~~~~~~~~~~~~~~~~~~~\tag 15
 $$
 Since $i_{(\dot {\phi}_t - \dot {\rho}_t)}\omega = du_t = i_{(X_{(u_t)})}\omega$ where $X_{u_t}$ 
 is the hamiltonnian vector field of $u_t$, we see that
 $$
 \dot {\phi}_t = \dot {\rho}_t + X_{u_t} = \dot {\rho}_t + (\rho_t)_* ((\rho_t)^{-1})_*( X_{u_t})
 $$
 Hence $\dot {\psi}_t = (\rho_t)^{-1})_*( X_{u_t}) = X_{(u_t \circ \rho_t)}$.
 This shows that $\psi_t$ is a hamiltonnian isotopy.
 
 \enddemo
  
 \pf
 
 \vskip .1in
 
 In formula (5), $\int_0^1 osc(u^{\Phi}_t) dt$ is nothing else than $l_H( \psi_t)$ and formula (5) can be written
 $$
 l(\Phi) = \int_0^1 |i(\dot {\rho}_t)\omega| dt + l_H( \psi_t))dt. ~~~~~~~~~~~~~~~\tag 5'
 $$
\vskip .2in 

{\bf 3. Prooof of theorem 1}

\vskip .1in

Clearly , $e(\phi) \geq 0$ for all $\phi$ and by definition $e(\phi) = e(\phi^{-1})$.

 To see that the triangular inequality holds, fix  a small positive number $\epsilon 
  \leq 1/8$ and a smooth
increasing  function $a: [0,1] \to [0,1]$ such that $a_{|[0,\epsilon)} = 0$ and a
$a_{| (1-\epsilon), 1]} = 1$ and let $\lambda(t) = a(2t)$ for $0 \leq t  \leq 1/2$
and $\mu(t) = a( 2t - 1)$ for $ 1/2 \leq t \leq 1$.

 \vskip .1in
 
 If $\Phi \in Iso(\phi)$ and $\Psi \in Iso(\psi)$ , we get an isotopy 
$\Phi * \Psi = (\sigma_t)\in Iso(\phi \psi)$ defined as:
$$
\sigma_t =\cases \phi_{\lambda(t)}, &\text{for $0 \leq t \leq 1/2$}\\
\phi_1\psi_{\mu(t)}, &\text{for $1/2\leq t \leq 1$}.\endcases
$$

\vskip .1in 

Let $c(\Phi, \Psi)$ be the set of all isotopies from $\phi\psi$ to the identity obtained as above.

 Clearly :

$$
e_0(\phi\psi)  \leq inf_{\Cal R} (l(\Cal R))
$$

where $\Cal R \in c(\Phi,\Psi)$. 

\vskip .1in

Since
$$
\dot {\sigma}_t  =\cases \lambda'\dot {\phi}_{\lambda(t)}, &\text{for $0 \leq t \leq 1/2$}\\
 \mu'\dot {\psi}_{\mu(t)}, &\text{for $1/2\leq y \leq 1$},\endcases
$$
\vskip .1in
we have:
$$
i(\dot {\sigma}_t)\omega = \cases \lambda'\Cal H^{\Phi}_{\lambda(t )}+ d(\lambda' 
u^{\Phi}_{\lambda(t)}) &\text{for $0 \leq t \leq 1/2$}\\
\mu'\Cal H^{\Phi'}_{\mu(t)} + d(\mu' u^{\Phi'}_{\mu(t)}),
 &\text{for $1/2\leq t \leq 1$},\endcases
$$

\vskip .1in

Therefore 

$$
l(\Phi*\Psi) = \int_0^{1/2} (|\lambda'\Cal H^{\Phi}_{\lambda(t )}|+ osc(\lambda' 
u^{\Phi}_{\lambda(t)})|)dt + \int_{1/2}^1(|\mu'\Cal H^{\Phi'}_{\mu(t)}| + osc(\mu' u^{\Phi'}_{\mu(t)})dt
$$
 By the change of variable formula , we get:
 
$$
l(\Phi*\Psi) = l(\Phi) + l(\Psi)
$$

Finally,

$$
e_0(\phi\psi)  \leq inf_{\Cal R} (l(\Cal R)) \leq inf_{\Phi} l(\Phi) + inf_{\Psi} l(\Psi) = e_0(\phi) + e_0(\psi).
$$
 
Therefore the triangular inequality holds true for $e_0$, and hence for $e$ as well. 

\vskip .1in

Showing that $e$ is non-degenerate is more delicate. Suppose that $e_0 (\phi) = 0$.

\vskip .1in

{\bf Step 1}

\vskip .1in

The statement $e_0(\phi) = inf (l(\Phi)) = 0$ means that for every $N$, there exists
 an isotopy $\Phi^N$ from $\phi$ to the identity such that $l(\Phi^N) \leq 1/N$.

Thus:
$$
\int_0^1 |\Cal H ^{\Phi^N}_t|dt \leq 1/N  ~~~~~~~~~~~~\tag 16
$$ 
and
$$
\int_0^1 osc (u^{\Phi^N}) dt \leq 1/N
$$
Hence
$$
 |\Cal H(\Phi^N)| = |\int_0^1 \Cal H ^{\Phi^N}_tdt| \leq  \int_0^1 |\Cal H ^{\Phi^N}_t|dt \leq   1/N.
$$

\vskip .1in

For any symplectic isotopy from $\phi$ to the identity $\Phi = (\phi_t)$, the 1-form
$$ 
\Cal H(\Phi) = \int_0^1  \Cal H ^{\Phi}_t dt
$$
is  the harmonic representative  of the cohomology class $\tilde {S}([\phi_t])$.
 
For any  symplectic isotopy  $\Phi =  (\phi_t)$  from $\phi$ to the identity
 $$
 \Cal H(\Phi^N) - \Cal H(\Phi) = \gamma(\Phi) \in \Gamma ~~~~~~~~~~~~~~~~~~~~~~~~~~~~\tag 17
 $$
 
since  $\Cal H(\Phi^N) - \Cal H(\Phi)$ is the harmonic representative of the image by $\tilde S$ of the
  class $[\phi^N_t * \phi_{(1-t)}]$ of the loop $\phi^N_t * \phi_{(1-t)}$.
 
 \vskip .1in

By (16) and (17), the distance $d(\Cal H(\Phi), \Gamma)$ from $\Cal H(\Phi)$ to $\Gamma$
satisfies:
$$
d(\Cal H(\Phi), \Gamma) \leq | \Cal H(\Phi) - (- \gamma(\Phi ^N)| = |\Cal H(\Phi^N)|\leq 1/N
$$
This says that $(\Cal H(\Phi))$ is arbitrarly close to $\Gamma$. Hence $(\Cal H(\Phi))\in \Gamma$.
This means that $\phi \in Ker S = Ham(M,\omega)$. 

\vskip .1in

The facts that $\Cal H (\Phi^N) \in \Gamma$ and $|\Cal H (\Phi^N)| \leq 1/N$ imply
that $\Cal H (\Phi^N) = 0$ for $N$ large enough since $\Gamma$ is discrete  (Ono's theorem).

Fix now  an isotopy $\Phi^N$ such that $\Cal H (\Phi^N) = 0$. To simplify the notations, we denote
 by $\Phi = (\phi_t)$ the isotopy $\Phi^N = (\phi^N_t)$. 

The Hodge decomposition of the isotopy $\phi_t$ gives:
$$
\phi_t = \rho_t \mu_t
$$
where $\rho_t$ is harmonic and $\mu_t$ is hamiltonian. We have:

$$
i(\dot {\phi}_t)\omega = \Cal H_t + du_t
$$
$$
\dot {\rho}_t =  (\Cal H_t)^\# = H_t,
$$
$$
 \int_0^1 \Cal H_t  dt = 0
 $$
 and
$$
 \int_0^1 (|\Cal H_t | + osc(u_t)) dt \leq 1/N 
$$

Hence
$$
\int_0^1 |\Cal H_t |dt \leq 1/N ;   \int_0^1  osc(u_t)dt \leq 1/N .~~~~~~~~~~~~~~~\tag 18 
$$

\vskip .1in

{\bf Step 2}

\vskip .1in

We are now going to deform the isotopy $\rho_t$ fixing the extremities  to a hamiltonian isotopy 
following [1], proposition II.3.1.

Let $Z_{(s,t)}$ be the family of symplectic vector fields:
$$
Z_{(s,t)} = t \dot {\rho}_{(s.t)} - 2s (\int_0^t (i(\dot {\rho}_u)\omega)du))\#. ~~~~~~~~~~~~~~\tag 19
$$
Clearly, $Z_{(0,t)} = 0$ and we have:
$$
\int_0^1 i(Z_{(s,t)})\omega ds = 0. ~~~~~~~~~~~~~~~~~~~~~~~~\tag 20
$$
Let $G_{(s,t)}$ be the 2-parameter family of diffeomorphisms obtained by integrating  $Z_{(s,t)}$
  with $t$ fixed, i.e. $G_{(s,t)}$  is  defined by the following equations:
  
 $$
 \frac {d}{ds}G_{(s,t)}(x) = Z_{(s,t)}(G^{-1}_{(s,t)}(x)) ,  G_{(0,t)}(x) = x.  ~~~~~~~~~~~~~\tag 21
 $$
 
 By (20), $G_{(1,t)}$ is a hamiltonnian diffeomorphism for all $t$. Since $Z_{(s,1)}  = \dot {\rho}_s - 2s
 ((\int_0^1 (i(\dot {\rho}_u)\omega)du) = \dot {\rho}_s$ , $s\mapsto G_{(s,1)}$ is an isotopy from the identity
  to $G{(1,1)} = \rho_1$. Hence the $g_t = G_{(1,t)}$ is an isotopy in $Ham (M,\omega)$ from $\rho_1$
 to the identity.

 \vskip .1in
 
 Consider the 2-parameter family of vector fields $V_{(s,t)}$ defined by:
 $$
 V_{(s,t)}(x) = \frac {d}{dt}G_{(s,t)}((G^{-1}_{(s,t)}(x))
 $$
 Clearly $\dot {g}_t =  V_{(1,t)}$.
 \vskip .1in
 We have ( see [1], proposition I.1.1):
 $$
 \frac {\partial}{\partial s}V_{(s,t)} = \frac {\partial}{\partial t} Z_{(s,t)}
  + [V_{(s,t)}, Z_{(s,t)}] ~~~~~~~~~~~~~~~~~~~~~~~~~~~~~~\tag 22
  $$

\vskip .1in

We will need the following 

\proclaim
{Proposition}

$$
i(V_{(1,t)})\omega = du_t
$$
where $u_t = \int_0^1 \omega(Z_{(s,t)}, V_{(s,t)}) ds$.

\endproclaim

\vskip .1in

\demo
{Proof}

>From equation 22
$$
0 = \frac {\partial}{\partial t} [\int_0^1 i(Z_{(s,t)})\omega ds] =
 \int_0^1 i(\frac {\partial}{\partial t}(Z_{(s,t)}))\omega ds 
$$
$$ 
 = \int_0^1 i(\frac {\partial}{\partial s}
 (V_{(s,t)}))\omega ds - \int_0^1 i([Z_{(s,t)}, V_{(s,t)}])\omega ds 
 = \int_0^1 (\frac {\partial}{\partial s} i(V_{(s,t)})\omega) ds-
  \int_0^1 i([Z_{(s,t)}, V_{(s,t)}])\omega ds
$$
$$
= i(V_{(1,t)})\omega - i(V_{(0,t)})\omega - \int_0^1 i([Z_{(s,t)}, V_{(s,t)}])\omega ds
= i(V_{(1,t)})\omega - d(\int_0^1 \omega(Z_{(s,t)}, V_{(s,t)}) ds)
$$
We used the facts that $V_{(0,t)} = 0$, $i([Z,V]\omega = L_Z i_V \omega - i_V L_Z\omega$ 
and $L_V\omega = L_V\omega = 0$.
\pf

\enddemo
 
 {\bf  Step 3: Norm estimates}

\vskip .1in

The harmonic vector fields $\dot {\rho}_t$ can be written as $\dot {\rho}_t = \sum_0^k \lambda_i(t)H_i$,
 where $H_i = h_i^{\#}$ and $(h_i)$ is a basis of harmonic 1-forms.
Formula (19) just says:
$$
Z_{(s,t)} =\sum_i  (t \lambda_i(st) -2s \int_0^t \lambda_i(u)du)H_i = \sum_i\mu_i(s,t)H_i ~~~~~~~~~~~~~~~~~~\tag 23
$$
Hence:
$$
|Z_{(s,t)}| = \sum_i|\mu_i(s,t)| \leq  t|\dot {\rho}_{st}| + 2s \int_0^t |\Cal H_t| dt
$$
$$
\leq t|\dot {\rho}_{st}| + 2s \int_0^1 |\Cal H_t| dt \leq t|\dot {\rho}_{st}| + 2s/N.
$$

\vskip .1in

On the other hand, we have:
$$
\omega(Z_{(s,t)}, V_{(s,t)}) = (i_{(Z_{(s,t)}}\omega)(V_{(s,t)}) = \sum_i \mu_i(s,t)h_i(V_{(s,t)})
$$

\vskip .1in

Consequently:
$$
|\omega(Z_{(s,t)}, V_{(s,t)})| \leq \sum_i |\mu_i(s,t) h_i(V_{(s,t)})|.
$$

Let $||h_i||$ be the sup norm of the 1-forms $h_i$, i.e $|| h_i|| = sup_{x \in M} |||h_i(x)|||$
and $|||h_i(x)|||$ is the norm of the linear map $h_i(x)$ on the tangent space $T_xM$.

We have:
$$
\sum_i |\mu_i(s,t)h_i(V_{(s,t)})|\leq (\sum_i |\mu_i(s,t)|) |V_{(s,t)}|E =  |Z_{(s,t)}||V_{(s,t)}|E $$

where $E = max \{||h_i||\}$. 

\vskip .1in
Hence
$$
|w_t| = |\int_0^1 \omega(Z_{(s,t)}, V_{(s,t)})ds| \leq \int_0^1 |\omega(Z_{(s,t)}, V_{(s,t)})|ds
$$
$$
\leq E(\int_0^1 (t|\dot {\rho}_{st}| + 2s/N)|V_{s,t})| ds.
$$

Let $A  = sup_{s,t} |V_{s,t}|$, then
$$
|w_t| \leq  AE\int_0^1 (t|\dot {\rho}_{st}| + 2s/N)ds = AE(\int_0^t (|\dot {\rho}_u| du) + 1/N) 
$$
$$
\leq AE(\int_0^1 (|\dot {\rho}_u| du + 1/N) \leq 2AE/N.
$$

Therefore $osc (w_t) \leq 4AE/N$, hence the length of the isotopy $\rho_t$ is
 less or equal to $4AE/N$,
and  therefore the Hofer norm of $\rho$ :  $||\rho||_H \leq  4AE/N$, where $\rho = \rho_1$.

\vskip .1in

{\bf Step 4}

\vskip .1in

Let $\Cal M$ denote the space of smooth maps $c:I = [0,1] \to W$, where $W$ is the space of
 symplectic vector fields on $(M,\omega)$ such that $c(0) = 0$ 
 with the Hofer norm
 $$
 ||c|| = \int_0^1 |c(t)| dt
 $$
 Here $|c(t)|$ is the norm given by formulas 13 and 14.
 
On the space $\Cal M \times I$ we define the distance $d(c,s),(c',s')) = ((||c-c'||^2 + (s-s')^2)^{1/2}$

Let $\Cal N$ be the space of smooth functions $u : I \times I \to U$, where $U$ is the space of symplectic 
vector fields with the metric $||u|| = sup_{s,t} |u(s,t)|$.

\vskip .1in

The family of vector fields  $V_{s,t}$ above is the image  of $\dot {\rho}_t$ by the following map:
$$
\Cal R : \Cal M \times I \to \Cal N
$$
where $\Cal R = \partial_t \circ I_s \circ a_s$ with

$a_s : c(t) \mapsto tc(st) - 2s ( \int_0^t i(c(u) \omega du)^\#$

$I_s : U_{s,t} \mapsto G_{s,t} : M \to M$ where the family of diffeomorphisms $G_{s,t}$ is obtained by 
integrating in $s$ like in formula 21.

and finally $\partial t : g_{s,t} \mapsto \partial/\partial t ( g_{s,t})$
 ( formula 22).

\vskip .1in

The mapping $\Cal R$ is a smooth map since all its components are smooth, consequently it is Lipschitz.
Therefore there is a constant $K$ 
such that $d(\Cal R(\dot {\rho}_t,s),(0,0))  = sup_{s,t} |V_{s,t}| \leq K (||\dot {\rho}_t||^2 + s^2)^{1/2}$
(Observe that $\Cal R(0,0) = 0$).

Therefore
$$
A = sup_{s,t} |V_{s,t}|\leq K ((1/N)^2 + s^2)^{1/2} \leq K ((1/N)^2 + 1)^{1/2}
\leq 2K
$$
Finally, we get:

$$
||\rho||_H \leq (4E(K ((1/N)^2 + 1)^{1/2}))/N
 \leq C/N
$$
where $C = 8EK$.

Remember now that $\phi = \rho \mu$ and $||\mu||_H \leq 1/N$. Therefore, $||\phi||_H \leq  (C +11)/N$
for all $N$. Hence $||\phi||_H = 0$ and consequently $\phi = id$.\pf

\vskip .1in

{\bf $Ham(M,\omega)$ is closed in $Symp(M,\omega)$}

\vskip .1in

Let $(h_n) \in Ham(M,\omega)$ be a sequence converging to $g\in Symp(M,\omega)$. There exists $N_0$ such 
that for all $N\geq N_0$, there exists an isotopy $\Phi^N \in Iso (g^{-1}h_N)$ with length 
$l(\Phi^N) \leq 1/N$. By step 1, $g^{-1}h_N$ is hamiltonian for N large. Hence $g$ is also hamiltonian.\pf

\vskip .1in

{\bf 4. Applications to the $C^0$ symplectic topology}

\vskip .1in

In [10], Oh and Muller defined the group of symplectic homeomorphisms, 
$Sympeo(M,\omega)$ as the closure of
the group $Symp(M,\omega)$ of 
$C^\infty$ symplectic diffeomorphisms of $(M,\omega)$ 
in the group $Homeo(M)$ of homeomorphisms of $M$ with the $C^0$ topology, and the group $Hameo(M,\omega)$ of
hamiltonian homeomorphisms. The group $Sympeo(M\omega)$ has only  $C^0$ topology induced from $Homeo(M)$, 
but $Hameo(M,\omega)$ has a more involved topology, called
 the {\it hamiltonian topology}, which combines 
the $C^0$ topology and the Hofer topology. 

\vskip .1in

Using our construction, we  define a {\it symplectic topology} 
on the space $Iso(M)$ of symplectic isotopies of $(M,\omega)$ as follows:

\noindent Fix a distance $d_0$ on $M$ (coming from some riemannian metric)
and define the distance $\overline {d}$ on the space $Homeo(M)$ of homeomorphismes of $M$ as
$$
\overline {d}(\phi,\psi) = max\{d(\phi,\psi), d(\phi^{-1},\psi^{-1})\}
$$
where 
$$
d(h,g) = max_x(d_0(h(x),g(x))
$$
for all $h,g \in Homeo(M)$.

 Then $(Homeo(M), \overline {d})$ is a complete
metric space and its metric topology is just the $C^0$ topology.
 On the space $\Cal PHomeo(M)$ of continuous paths
$\lambda :[0,1] \to Homeo(M)$, we put the metric topology from the distance 
$$
\overline {d}(\lambda,\mu) = sup_{t\in [0,1]}\overline {d}(\lambda(t),\mu(t)).
$$
 We define the {\it symplectic distance } on $Iso(M)$ by:
$$
 d_{symp}(\Phi, \Psi) = \overline {d}(\Phi,\Psi) + D(\Phi.\Psi)
$$
where $D$ is given by formula (5').

\noindent We call the {\it symplectic topology} on $Iso(M)$ the metric topology 
defined by the above distance. This topology reduces to the "hamiltonian topology" of [10]
on paths in $Ham(M,\omega)$.

\vskip .1in

We now define a set (which we 
conjecture to be a group)
$SSympeo(M,\omega)$ as follows: $h \in SSympeo(M)$ iff
   there exists a continuous path
   $\lambda : [0,1] \to Homeo(M)$ such that $\lambda(0) = id ;  \lambda(1) = h$
   and a sequence $\phi_n^t$ of symplectic isotopies, which converge  
   to $\lambda$ in the
   $C^0$ topology ( induced by the norm $\overline {d}$) and such that $D((\Phi_j (\Phi_i)^{-1})$ tends
   to zero when $i$ and $j$ go to infinity. This set, which contains $Hameo(M,\omega)$
    will play a major role in the $C^0$
symplectic topology, where one may consider the subgroup of $Homeo(M)$
    it generates. This subgroup contains $Hameo(M,\omega)$ as a normal subgroup.
 \vskip .2in   
    
{\bf 5. Final Remarks}

\vskip .1in

The metric $e$ obtained here is not an " extension" of the Hofer
 metric since we do not know if $e(\phi) = ||\phi||_H$ when $\phi \in Ham(M,\omega)$. We only know
 that $e(\phi) \leq ||\phi||_H$. The problem of extending the Hofer norm was considered in [3]. Here 
 we would like to make some remarks about the results of [3].
 
 \vskip .1in

{\bf Extension of Oh's spectral norm}. It is obvious that formulas
 of the extensions of the Hofer metric given  in [3]  give 
in fact extensions for any bi-invariant metric  on $Ham(M,\omega)$.
 Theorem 2 in [3] uses only the properties of bi-invariance and not
  the Hofer norm. Then theorem 2 of [3] can be rephrased as
 
 \vskip .1in
 
 \proclaim
 {Theorem 2}
 
 Let $(M,\omega)$ be a symplectic manifold such that the homomorphism $S$  
  admits a continuous homomorphic right inverse,
 then any bi-invariant metric on $Ham(M,\omega)$ extends to a right invariant metric on $Symp(M,\omega)$.
 
\endproclaim

\vskip .1in

Under the hypothesis of the theorem above, the spectral
norm $||.||_{\Cal O}$ of Oh extends to all of $Symp(M,\omega)_0$. For the definition of Oh's spectral
 norm, we refer to [9].
An example where this  hypothesis holds is $T^{2n}$ with its natural symplectic form.

\vskip .in

\proclaim
{Theorem 3}

If $\Gamma = 0$, Oh's spectral distance extends to $Symp(M,\omega)_0$.

\endproclaim

\vskip .1in

\demo
{Proof}

Let $\phi_i, i= 1,2$ two symplectomorphisms and $\Phi_i =(\phi^i_t) 
\in Iso (\phi_i)$. The harmonic 1-forms $\Cal H(\Phi_i)$ depend only on $\phi_i$.
 Let $\rho_i$ be the time one of the 1-parameter group generated by $\Cal H(\Phi_i)^\#$, 
 then $\psi_i = \phi_i \rho_i^{-1} \in Ham(M,\omega)$. We define the Oh distance $d_{\Cal O}$ of $\phi_1$
  and $\phi_2$ by:
 $$ 
d_{\Cal O}(\phi_1, \phi_2) = |\Cal H(\Phi_1) - \Cal H(\Phi_2)| + ||\psi_1 \psi_2 ^{-1}||_{\Cal O}.
$$

The cases where $\Gamma = 0$  include oriented compact surfaces of genus bigger than one. More recently,
 Kedra, Kotschick 
and Morita [6] found a longer list
of compact symplectic manifolds with vanishing flux group. 

\vskip .1in

 \noindent{\bf References}
\parindent=0pt   

[1]~ A. Banyaga {\it ~Sur la structure du groupe des diff\'eomorphismes qui
pr\'eservent une forme symplectique},
Comment. Math. Helv. 53(1978) pp.174--227.

[2]~ A. Banyaga {\it ~The structure of classical diffeomorphisms groups},
Mathematics and its applications vol 400.
Kluwer Academic Publisher's Group, Dordrecht, The Netherlands (1997).

[3]~ A. Banyaga , P. Donato {\it ~ Lengths of Contact Isotopies and Extensions of the Hofer Metric}
Annals of Global Analysis and Geometry 30(2006) 299-312

[4]~ Z. Han {\it Bi-invariant metrics on the group of symplectomorphisms}, Preprint

[5]~H. Hofer {\it~On the topological properties of symplectic maps}, Proc.
Royal Soc. Edimburgh 115A (1990), pp.25--38

[6] J. Kedra, D. Kotschick, S. Morita {\it Crossed flux homomorphism and vanishing theorem for flux groups}
Geom Funct. Anal 16(2006)no 6 1246-1273.

[7]~F. Lalonde, D. McDuff {\it~ The geometry of symplectic energy}, Ann.
Math. 141 (1995) 349 - 37 .

[8]~F. Lalonde, L.Polterovich {\it Symplectic diffeomorphisms as isometries of Hofer's norm}
Topology 36(1997) 711-727

[9]~Y-G. Oh {\it Spectral invariants, analysis of the Floer moduli space,
 and the geometry of hamiltonian diffeomorphisms}, Duke Math. J. 130(2005) 199-295

[10]~Y-G. Oh and S. Muller {\it The group of hamiltonian homeomorphisms and $C^0$-symplectic topology}
J. Symp.Geom. to appear
 
[11] ~K. Ono {\it Floer-Novikiv cohomology and the flux conjecture}
Geom. Funct. Anal.16(2006) no 5 981-1020
 
[12]~L. Polterovich {\it~ Symplectic displacement energy for Lagrangian
submanifolds}, Erg. Th and Dyamical Systems 13
(1993), 357-367.

[13]~C. Viterbo {\it~ Symplectic topology as the geometry of generating
functions}, Math. Annalen 292(1992), 685-710

[14]~ F. Warner {\it Foundations of differentiable manifolds and Lie groups}
Scott, Foresman and Company (1971).

 \vskip .3in
\baselineskip 12pt
\noindent Department of Mathematics \newline
\noindent The Pennsylvania State University \newline
\noindent University Park, PA 16802

\end